\documentclass[11pt,a4paper,reqno]{amsart}

\usepackage{mathrsfs,amssymb}

\newtheorem{theorem}{Theorem}
\newtheorem{lemma}[theorem]{Lemma}

\newtheorem{corollary}[theorem]{Corollary}

\newtheorem{question}[theorem]{Question}

\theoremstyle{definition}

\newtheorem{definition}[theorem]{Definition}

\theoremstyle{plain}

\newcommand{\BZ}{\mathbb{Z}}

\newcommand{\ol}[1]{\overline{#1}} 

\begin{document}
    \title[Word problems and monoid automata]{Word problems recognisable by deterministic blind monoid automata}

\maketitle

\begin{center}
    Mark Kambites \\

    \medskip

    Fachbereich Mathematik / Informatik, \  Universit\"at Kassel \\ 
    34109 Kassel, \  Germany

    \medskip

    \texttt{kambites@theory.informatik.uni-kassel.de} \\
    
\end{center}   

\begin{abstract}
We consider blind, deterministic, finite automata equipped with a register
which stores an element of a given monoid, and which is modified by right
multiplication by monoid elements. We show that, for monoids $M$ drawn from
a large class including groups, such an automaton accepts the word problem
of a group $H$ if and only if $H$ has a finite index subgroup which embeds
in the group of units of $M$. In the case that $M$ is a group, this answers a question of Elston
and Ostheimer.
\end{abstract}

\section{Introduction}

Several authors have studied finite automata
augmented with a memory register, which stores at any moment an element
of a given monoid $M$. Such an automaton, which we shall term an \textit{$M$-automaton}, cannot read its register, but
it may change the contents of the register by multiplying on the right by
some element of the monoid. The register is initialised with the identity
element, and the automaton accepts an input word $w$ exactly if, by
reading this word, it can reach a final state, in which the register
has returned to the identity element.

A number of notable classes of languages can be characterised in this way.
If $M$ is taken to be a free abelian group of rank $n$, then the languages
accepted by $M$-automata are exactly those accepted by the blind $n$-counter
machines studied
by Greibach \cite{Greibach78}. As another example, Gilman \cite{Gilman91}
has observed that if $M$ is a polycyclic monoid of rank $2$ or more, then we
obtain exactly the context-free languages. Other examples have been studied,
implicitly or explicitly, by Ibarra, Sahni and Kim \cite{Ibarra76}, by Dassow
and Mitrana \cite{Dassow00} and by Mitrana and Stiebe \cite{Mitrana97,Mitrana01}.

At the same time, an area of perennial interest in combinatorial and
computational algebra is the study of \textit{word problems} of groups.
It is well-known that many structural properties of groups are reflected
in the language theoretic properties of their word problems, and vice
versa. It seems very
natural to ask if there is any connection between the structural properties
of a given monoid $M$ and of the collection of groups whose word problems
are recognised by $M$-automata.

This problem was explicitly considered by Elston and Ostheimer \cite{Elston04},
in the case where the register monoid is a group. They showed that a group $H$
has word problem recognisable by a deterministic $G$-automaton with a
certain inverse property if and only $H$ has a finite index subgroup which
embeds in $G$. They also considered two possible ways in which one might think to strengthen this
theorem. Firstly, they demonstrated that the hypothesis of determinism cannot
in general be removed, by exhibiting a non-deterministic
$(F_3 \times F_3)$-automaton accepting the word problem of the free
abelian group $\BZ^3$.

Secondly, they posed the question of whether the inverse property
hypothesis can be removed, that is, whether a group $H$ whose word problem
is accepted by a deterministic $G$-automaton necessarily has a finite index
subgroup which embeds in $G$. The primary objective of this note is to provide
a positive answer to this question. A secondary intention is to show, in addition,
that one need not assume the register monoid to be a group; rather, it
suffices that it belongs to a large class of monoids satisfying a weak right
cancellativity condition. Specifically, we prove the following (definitions
can be found in Section~\ref{sec:basics} below).
\begin{theorem}\label{thm:maintheorem}
Let $M$ be a monoid with unique left inverses and $H$ a finitely generated
group. Then the word problem for $H$ is accepted by a deterministic
$M$-automaton if and only if $H$ has a finite index subgroup which embeds
in the group of units of $M$.
\end{theorem}

In addition to this introduction, this paper comprises three sections.
In Section~\ref{sec:basics} we recall some basic definitions and results,
while Section~\ref{sec:mainproof} is devoted to the proof of
Theorem~\ref{thm:maintheorem}. Finally, in Section~\ref{sec:consequences}
we discuss some consequences of our result, and pose a related
question which remains open.

\section{Basics}\label{sec:basics}

In this section we recall the basic definitions which are used in
the sections that follow. We assume a basic familiarity with finite
automata and with group theory; the reader lacking prior experience in
either of these areas is advised to consult one of the many textbooks
on each area, such as \cite{Lawson03} or \cite{Robinson96}.

Let $M$ be a monoid with identity $1$, and $X$ be a finite alphabet. By an
\textit{$M$-automaton $A$ over $X$} , we mean a finite automaton $A$ over
$M \times X^*$. We say that $A$ \textit{accepts} a word $w \in X^*$
if it accepts $(1,w) \in M \times X^*$ when considered as a 
finite automaton in the usual sense, that is, if there is a path from the
start state to a terminal state labelled $(1,w)$. The \textit{language
recognised} or \textit{accepted} by $A$ is the set of all words $w \in X^*$
which are accepted by $A$.

The \textit{underlying automaton} of $A$ is the finite automaton over $X$
which is obtained from $A$ by disregarding the $M$ component of the edge
labels. It is easily seen to accept a regular language which contains the
language accepted by $A$.

We say that an $M$-automaton $A$ is \textit{deterministic} if its edges are
labelled by elements of $M \times X$, and for every state $q$ and letter
$x \in X$ there is at most one edge leaving $x$ with a label of the form
$(g,x)$ for some $g \in M$.

Of particular interest to us is the case where the alphabet $X$ is viewed
as a monoid generating set for a finitely generated free group $F$. For
a word $w \in X^*$ we denote by $\overline{w}$ the element of $F$ represented.
We say that a $G$-automaton $X$ over $F$ \textit{accepts a subset
$Q \subseteq F$} if it accepts the preimage of $Q$ in $X^*$, that is, the
language of \textbf{all} words in $X^*$ which represent elements of $Q$.

If $H$ is a group generated as a monoid by a subset $X$, then $H$ is
naturally isomorphic to $F / Q$ for some free group $F$ generated as a
monoid by $X$ and some normal subgroup $Q$ of $F$. The \textit{word
problem} for $H$ is the set of all words in the free monoid $X^*$ which
represent the identity in $H$, that is, which represent elements of $Q$ in $F$.

Let $M$ be a monoid with identity $1$. Recall that the \textit{group of
units} $G(M)$ of $M$ is the set of elements in $M$ which have a two-sided
inverse with respect to $1$, that is, the largest subgroup of $M$ containing
$1$.

We say that $M$ has \textit{unique left inverses} if whenever
$a,b,c \in M$ are such that $ba = 1 = ca$, we have $b = c$.
This condition can be regarded as an extremely weak cancellativity 
property. In particular, it is satisfied by all right cancellative
monoids, and hence by all groups. It is also satisfied by a number of
other monoids of particular interest in this context; these include
the polycyclic monoids, which will be discussed further in
Section~\ref{sec:consequences}.

\section{Proof of the Main Theorem}\label{sec:mainproof}

In this section we prove our main theorem. The proof requires a number of 
definitions and lemmas. Throughout this section, we assume that $M$ is
a monoid with identity $1$; from a certain point on, we shall require that
it has unique left inverses. The reader better versed in group theory than semigroup theory
may prefer, at a first reading, to imagine that $M$ is a group throughout.
We begin with a straightforward definition.
\begin{definition}\label{def_tsm}
An [deterministic] $M$-automaton $A$ accepting a language $L$ is called
[deterministic] \textit{terminal state minimal} or \textit{[D]TSM} if $L$
is accepted by no [deterministic] $M$-automaton with strictly fewer 
terminal states.
\end{definition}
Since the natural numbers are well-ordered, it is immediate that any
language accepted by a [deterministic] $M$-automaton is accepted by a
[D]TSM $M$-automaton.

\begin{lemma}\label{lem_terminalsused}
Let $A$ be a TSM or DTSM $M$-automaton accepting a language $L$. Let $q_0$
be the initial state of $A$ and let $q$ be a terminal state of $A$. Then
there is a path from $q_0$ to $q$ with label $(1, w)$ for some $w \in L$.
\end{lemma}
\begin{proof}
If not, then the automaton obtained from $A$ by making $q$ non-terminal
accepts the same language with strictly fewer terminal states, contradicting
the assumption that $A$ is TSM or DTSM.
\end{proof}

Now let $F$ be a free group generated as a monoid by a finite subset $X$.
Let $Q$ be a subgroup of $F$, and suppose that $Q$ is accepted by a
TSM or DTSM $M$-automaton $A$.

\begin{lemma}\label{lem_changeinitial}
Let $q$ be a terminal state of $A$, 
and let $A_q$ be the automaton which is the same as $A$, except that the
initial state is $q$. Then $A_q$ accepts a subset of the words representing
elements of $Q$. Moreover, if $A$ is deterministic then $A_q$ accepts
exactly the words representing elements of $Q$.
\end{lemma}
\begin{proof}
Let $q_0$ be the initial state of $A$. By Lemma~\ref{lem_terminalsused},
there is a path labelled $(1, w)$ from $q_0$ to $q$, for some $w \in X^*$
such that $\ol{w} \in Q$. 
Now suppose $x \in X^*$ is accepted by $A_q$. Then there is a path from $q$ 
to a terminal state labelled $(1,x)$. It follows that there is a path from 
$q_0$ to the same terminal state labelled $(1, wx)$, so that $\ol{wx} \in Q$. 
But since $Q$ is a subgroup and $\ol{w} \in Q$, it follows that $\ol{x} \in Q$.

For the converse, we must suppose that $A$ is deterministic. Now if $x \in X^*$
is such that $\ol{x} \in Q$ then $\ol{wx} \in Q$. So $wx$ is accepted by $A$,
that is, there is a path from $q_0$ to a terminal state labelled $(1, wx)$. But since the automaton is
deterministic, this path must reach $q$ after reading $w$, so there is a path
from $q$ to a terminal state labelled $(1, x)$. But this path also exists in
$A_q$, so $x$ is accepted by $A_q$.
\end{proof}

\begin{corollary}\label{cor_betweenterminals}
Suppose $A$ is deterministic. Then the language of all words $w \in X^*$ such
that $(1, w)$ labels a path between two terminal states is exactly the
membership language of the subgroup $Q$.
\end{corollary}
\begin{proof}
The empty word over $X$ represents the identity in $F$, which must be a
member of the subgroup $Q$. Since the automaton is deterministic, it
follows that the initial state is a terminal state. Hence, the given
language contains all words $w \in F$ such that $(1, w)$ labels a path from the initial state to a terminal
state, that is, the membership language of $Q$. Conversely, if $(1,w)$ labels
a path between terminal states $p$ and $q$ then $w$ is accepted by $A_p$, so
by Lemma~\ref{lem_changeinitial}, $\ol{w} \in Q$ as required.
\end{proof}

\begin{corollary}\label{cor_terminalsconnected}
Suppose $A$ is deterministic, and let $p$ and $q$ be terminal states
of $A$. Then there exists $w \in X^*$ with $\ol{w} \in Q$ such that
there is a path from $p$ to $q$ labelled $(1, w)$.
\end{corollary}
\begin{proof}
By Lemma~\ref{lem_changeinitial}, the automaton $A_p$ with initial state $p$ accepts
$Q$. Clearly $A_p$ is still DTSM, so by Lemma~\ref{lem_terminalsused},
there is a path in $A_p$ from $p$ to $q$ labelled $(1, w)$ for some $w \in X^*$
with $\ol{w} \in Q$. But this path also exists in $A$.
\end{proof}

\begin{lemma} \label{lem_welldef}
Suppose $M$ has unique left inverses, $A$ is deterministic and $Q$ is
a normal subgroup of $F$. Suppose $(g, w)$ and $(h, z)$ both
label paths between
(possibly different) pairs of terminal states. Suppose further that
$\ol{w} Q = \ol{z} Q$. Then $g = h$.
\end{lemma}
\begin{proof}
Suppose $(g,w)$ labels a path from $a$ to $b$, while $(h, z)$ labels a
path from $c$ to $d$. By Lemma~\ref{lem_changeinitial}, we can assume
without loss of generality that $a$ is the initial state.

By Corollary~\ref{cor_terminalsconnected}, there is a path $x$ labelled
$(1, u)$ from $a$ to $c$, and a path labelled $(1, v)$ from $d$ to $b$, 
for some $u, v \in X^*$ such that $\ol{u}, \ol{v} \in Q$. Hence, there is
a path $y$ labelled $(1,u) (h,z)  (1,v) = (h, uzv)$ from $a$ to $b$. Now
since $\ol{u}, \ol{v} \in Q$ and $Q$ is normal, we must have
that $(\ol{uzv})(\ol{z})^{-1} \in Q$. Let $t \in X^*$ be a word representing
$(\ol{z})^{-1}$. Then the automaton must accept $u z v t$. Also,
by assumption, we have $\ol{w} Q = \ol{z} Q$, so that
$(\ol{w}) (\ol{z})^{-1} = \ol{w t} \in Q$,
so the automaton must accept
$w t$. Hence, there are paths from $a$ to some terminal
state, labelled $(1, uzvt)$ and $(1, wt)$. Since the
automaton is deterministic, these paths must begin by following the
paths $x$ and $y$ respectively, until they reach state $b$. Now there
must be paths from $b$ to a terminal state labelled $(g', t)$ and
$(h', t)$, where $g'$ and $h'$ are right inverses of $g$ and $h$
respectively. But the automaton is deterministic, so it follows that 
$g' = h'$, and now by the unique left inverse condition, that $g = h$ as 
required.
\end{proof}

We assume from here on that $A$ is deterministic, $M$ has unique
left inverses, and $Q$ is normal.  Now let $J$ be the set of all words
in $X^*$ which are read between terminal states in the underlying finite
state automaton of $A$. Let
$K = \lbrace j Q \mid j \in J \rbrace \subseteq F / Q$.

\begin{lemma} \label{lemma_ksubgroup}
$K$ is a subgroup of $F / Q$, with index bounded above by the number
of states in $A$.
\end{lemma}
\begin{proof}
Suppose $j, k \in J$ so that $\ol{j} Q, \ol{k} Q \in K$. Suppose $j$ labels
a path between terminal states $a$ and $b$ in the underlying automaton,
while $k$ labels a path between terminal states $c$ and $d$. By
Corollary~\ref{cor_terminalsconnected}, there is a path from $b$ to $c$
labelled by some word $q \in X^*$ such that $\ol{q} \in Q$. Hence, there
is a path from $a$ to $d$ labelled $jqk$, so that $j q k \in J$ and
$\ol{jqk} Q \in K$. But
clearly $\ol{j q k} Q = \ol{j} Q \ol{k} Q$, so that $\ol{j} Q \ol{k} Q \in K$.
Thus, $K$ is closed under composition.

Now suppose $j \in J$ so that $\ol{j} Q \in K$. Then there is a path in the
underlying automaton between terminal states labelled $j$, so there is
a path $p$ in $A$ between terminal states $a$ and $b$ labelled $(g, j)$ for
some $g \in M$. By Lemma~\ref{lem_changeinitial}, we can assume without
loss of generality that $a$ is the initial state of $A$. Let $t \in X^*$
be a word representing $(\ol{j})^{-1}$. Now certainly $\ol{j t} = 1 \in Q$,
so there must be a path from $a$ to a terminal state labelled $(1, j t)$.
Since $A$ is deterministic, this path must follow $p$ as far as state $b$,
so there must be a path from $b$ to a terminal state labelled $(g', t)$,
where $g'$ is a right inverse of $g$ in $M$.
Hence, there is a path between terminal states in the underlying automaton,
labelled $t$. Thus, $t \in J$ and so $(\ol{j})^{-1} Q = \ol{t} Q \in K$.

We have shown that $K$ is a subgroup of $F / Q$. It remains to justify the
claim about its index. Recall that a state $q$ of $A$ is \textit{accessible}
if there is a path in $A$ from the initial state to $q$. For each accessible
state $q$ in $A$, choose $g_q \in M$ and
$w_q \in X^*$ such that $(g_q, w_q)$ labels a path from the initial state
of $A$ to state $q$. For each word $w_q$, let $v_q \in X^*$ be a word
representing $(\ol{w_q})^{-1}$. Now $\ol{w_q} \ol{v_q} = 1 \in Q$,
so that $(1, w_q v_q)$ labels a state from the initial state to a
terminal state. Since $A$ is deterministic, this must pass through $q$,
and so there must be a path from $q$ to a terminal state labelled
$(g_q', v_q)$ for some right inverse $g_q'$ of $g_q$.

Now let $v \in F$. Then there is a path in $A$ labelled $(g, v)$ for some
$g \in M$, from the initial state to some state $q$. Now $q$ is clearly
accessible, so there is a path from $q$ to a terminal state labelled $(g_q', v_q)$.
It follows that the underlying finite automaton of $A$ accepts the word
$(v v_q)$, so that $v v_q \in J$. But now
$(\ol{v}) (\ol{w_q})^{-1} Q = \ol{v} \ol{v_q} Q \in K$, so that $\ol{v} Q$
and $\ol{w_q} Q$ lie in the same coset of $K$. Hence, the number of cosets of
$K$ in $F / Q$ is bounded above by the number of states in $A$, as required.
\end{proof}

We are now in a position to define the desired embedding of $K$ into
the group $G(M)$.

\begin{lemma}\label{lemma_embedding}
There is a well-defined embedding $\sigma : K \to G(M)$ given by
setting $\sigma(\ol{w} Q)$ to be the unique $g \in M$ such that
$(g, w)$ labels a path between terminal states in $A$.
\end{lemma}
\begin{proof}
That $\sigma$ is well-defined function from $K$ to $M$ follows directly
from Lemma~\ref{lem_welldef}.

To see that $\sigma$ is a homomorphism, suppose $j, k \in J$.
As in the previous proof, suppose $(g, j)$ labels
a path between terminal states $a$ and $b$ in $A$, while
while $(h, k)$ labels a path between terminal states $c$ and $d$. By
Corollary~\ref{cor_terminalsconnected}, there is a path from $b$ to
$c$ labelled by $(1, q)$ for some element $q \in X^*$ with $\ol{q} \in Q$.
Hence, there is 
a path from $a$ to $d$ labelled $(gh, jqk)$. It follows that
$\sigma(\ol{j} Q) \sigma(\ol{k} Q) = gh = \sigma(\ol{jqk} Q) = \sigma(\ol{j} Q \ol{k} Q)$,
as required.

Now if $j \in X^*$ is such that $\ol{j} Q = Q$ is the identity in $K$,
then we have $\ol{j} \in Q$. By the original assumption on the
automaton, there is a path between terminal states labelled $(1,j)$, so
that $\sigma(Q) = \sigma(\ol{j} Q) = 1$ in $M$.  Since the class of groups
is closed under the taking of semigroup homomorphic images, the image of
$H$ under $\sigma$ is a subgroup of $M$. Since it contains the identity of
$M$ is must be contained in the group of units $G(M)$.

It remains to show that $\sigma$ is injective. Suppose $j \in J$ is such that
$\sigma(\ol{j} Q) = 1$. By the definition of $\sigma$, there is a path between
terminal states $a$ and $b$ with label $(1, j)$. By Lemma~\ref{lem_changeinitial},
we may assume without loss of generality that $a$ is the initial state. It
follows that $j$ is accepted by the $G$-automaton $A$, so by our original
assumptions, $\ol{j} \in Q$ and $\ol{j} Q$ is the identity in $K$. Thus,
$\sigma$ is injective.
\end{proof}

The above results combine with a result from \cite{Elston04} to prove our
main theorem.
\begin{theorem}
Let $M$ be a monoid with unique left inverses and $H$ a finitely generated
group. Then the word problem for $H$ is accepted by a deterministic
$M$-automaton if and only if $H$ has a finite index subgroup which embeds
in the group of units of $M$.
\end{theorem}
\begin{proof}
If $L$ is finitely generated group $H$ with word problem accepted by an
$M$-automaton, then it follows from Lemmas~\ref{lemma_ksubgroup} and
\ref{lemma_embedding} that $H$ has a finite index subgroup $K$ which
embeds in the group of units of $M$.

Conversely, if $H$ has a finite index subgroup which embeds in the group
of units $G(M)$ of $M$, then by \cite[Theorem~7]{Elston04}, the word problem
of $H$ is accepted by a deterministic $G(M)$-automaton, and hence by a
deterministic $M$-automaton.
\end{proof}

\section{Consequences and Questions}\label{sec:consequences}

An immediate corollary of our main theorem is that if $M$ is a monoid
with unique left inverses and finite group of units, then no deterministic
$M$-automaton accepts the word problem of an infinite group.

One example is of particular interest. Let $X$ be an alphabet of size
$n$. Recall that the \textit{polycyclic monoid} of rank $n$ is the submonoid of
the semigroup of partial bijections on $X^*$ generated by the functions of the form
$$P_x : X^* \to X^*x, \ w \mapsto wx$$
where $x \in X$ and their relational inverses
$$Q_x : X^*x \to X^*, \ wx \mapsto w.$$
It is the natural semigroup-theoretic model of the operations which can
be performed on a pushdown stack with alphabet $X$; the element $P_a$
corresponds to ``pushing'' $a$, and $Q_a$ to ``popping'' $a$ from the
stack. More details of this viewpoint are contained in
\cite[Section~7]{Gilman91}, while a more general discussion of
polycyclic monoids can be found in \cite[Chapter~9]{Lawson98}.

It is readily verified that no non-identity element of a polycyclic
monoid has both a left and a right inverse, so that the group of units
is trivial. One can
also check that such a monoid has the unique left inverse property. It
follows, then, that no word problem of an infinite group is accepted by
a deterministic polycyclic monoid automaton. This contrasts with the
situation for non-deterministic polycyclic monoid automata. We have already
noted that the latter recognise exactly the context-free languages \cite{Gilman91};
it follows by a well-known theorem of Muller and Schupp \cite{Muller83}, that
they can recognise word problems of exactly the finitely generated virtually
free groups.

We mentioned above that Elston and Ostheimer \cite{Elston04} have exhibited
an example of non-deterministic $(F_3 \times F_3)$-automaton which recognises
the word problem of $\BZ^3$. Since $\BZ^3$ does not have a finite index
subgroup which embeds in $F_3 \times F_3$, it follows by Theorem~\ref{thm:maintheorem} (or by
\cite[Theorem~7]{Elston04}) that even when the register monoid is a group
$G$, deterministic and non-deterministic $G$-automata do not in general
accept the same word problems.

However, a consequence of Theorem~\ref{thm:maintheorem}, together with
the Muller-Schupp theorem and \cite[Lemma~3.2]{Dassow00}, is that a
group word problem accepted by a free group automaton is always accepted
by a deterministic free group automaton. One is drawn to ask the following
question.
\begin{question}
For what monoids $M$ is it true that deterministic and non-deterministic
$M$-automata can accept the same group word problems?
\end{question}

\section*{Acknowledgements}

This research was supported by a Marie Curie Intra-European Fellowship
within the 6th European Community Framework Programme. The author would like
to thank John Fountain and Gretchen Ostheimer for some helpful conversations,
and Kirsty for all her support and encouragement.

\bibliographystyle{plain}

\end{document}